\documentclass[12pt,reqno]{amsart}

\usepackage{setspace}

\usepackage{enumitem} %this enables references to correct labels in \enumerate

\usepackage
{hyperref, epigraph}

\usepackage{breakurl}   %this allows for url breaks in pdf

\theoremstyle{definition}

\DeclareMathOperator{\N}{{\mathbb N}}

\DeclareMathOperator{\Q}{{\mathbb Q}}
\DeclareMathOperator{\R}{{\mathbb R}}

\numberwithin{equation}{section}

\newcommand\astr{{{}^\ast\hspace{-2.5pt}\R}}

\newcommand{\hr} {{{}^{\mathfrak{h}}\hspace*{-2.3pt}\R}}

\newcommand{\st}{\textbf{st}}

\author[J. Bair]{Jacques Bair}\address{J. Bair, HEC-ULG, University of
Liege, 4000 Belgium}\email{j.bair@ulg.ac.be}

\author[P. Blaszczyk]{Piotr B\l{}aszczyk}\address{P. B\l{}aszczyk, Institute
of Mathematics, Pedagogical University of Cracow,
Poland}\email{pb@up.krakow.pl}

\author[P. Heinig]{Peter Heinig} \address{P. Heinig}
\email{peter.c.heinig@gmail.com}

\author[V. Kanovei]{Vladimir Kanovei} \address{V. Kanovei, IPPI RAS,
Moscow, Russia}\email{kanovei@googlemail.com}

\author[M. Katz]{Mikhail G. Katz}\address{M. Katz, Department of
Mathematics, Bar Ilan University, Ramat Gan 52900
Israel}\email{katzmik@macs.biu.ac.il}

\begin{document}

%\doublespacing

\thispagestyle{empty}

%\huge

\title{19th century real analysis, forward and backward}

\begin{abstract}
19th century real analysis received a major impetus from Cauchy's
work.  Cauchy mentions variable quantities, limits, and
infinitesimals, but the meaning he attached to these terms is not
identical to their modern meaning.  

Some Cauchy historians work in a conceptual scheme dominated by an
assumption of a teleological nature of the evolution of real analysis
toward a preordained outcome.  Thus, Gilain and Siegmund-Schultze
assume that references to \emph{limite} in Cauchy's work necessarily
imply that Cauchy was working with an Archi-medean continuum, whereas
infinitesimals were merely a convenient figure of speech, for which
Cauchy had in mind a complete justification in terms of Archimedean
limits.  However, there is another formalisation of Cauchy's
procedures exploiting his \emph{limite}, more consistent with Cauchy's
ubiquitous use of infinitesimals, in terms of the \emph{standard part
principle} of modern infinitesimal analysis.  

We challenge a misconception according to which Cauchy was allegedly
forced to teach infinitesimals at the \emph{Ecole Polytechnique}.  We
show that the debate there concerned mainly the issue of \emph{rigor},
a separate one from \emph{infinitesimals}.  A critique of Cauchy's
approach by his contemporary de Prony sheds light on the meaning of
rigor to Cauchy and his contemporaries.  An attentive reading of
Cauchy's work challenges received views on Cauchy's role in the
history of analysis, and indicates that he was a pioneer of
infinitesimal techniques as much as a harbinger of the
\emph{Epsilontik}.

Keywords: butterfly model; continuity; infinitesimals; \emph{limite};
standard part; variable quantity; Cauchy; de Prony
\end{abstract}

\maketitle

\tableofcontents

\epigraph{Since\, Weierstrass's\, time, we\, have held a fairly
contemptuous view of the infinitesimalists which I regard as unfair.
--\;Ivor Grattan-Guinness}

\section{Introduction}
\label{s1b}

Cauchy exploited the concepts of variable quantity, limit, and
infinitesimal in his seminal 1821 textbook \emph{Cours d'Analyse}
(CdA).  However, the meaning he attached to those terms is not
identical to their modern meanings.  While Cauchy frequently used
infinitesimals in CdA, some scholars have argued that Cauchyan
infinitesimals are merely shorthand for prototypes of
$\epsilon,\delta$ techniques.  Moreover, one can legitimately ask
whether the material found in CdA was actually taught by Cauchy in the
classroom of the \emph{Ecole Polytechnique} (EP).  A valuable resource
that sheds information on such issues is the archive of summaries of
courses and various \emph{Conseil} meetings at the EP, explored by
Guitard (\cite{Gu86}, 1986), Gilain (\cite{Gi89}, 1989), and others.
Among the key figures at EP at the time was Gaspard de Prony, whose
critique of Cauchy's teaching will be examined in Sections~\ref{s37}
and \ref{s311}.  While de Prony was critical of Cauchy, a careful
examination of the criticism indicates that de Prony's main target was
what he felt was excesssive rigor, rather than an alleged absence of
infinitesimals.  While scholars sometimes claim that Cauchy avoided
infinitesimals in the 1820s, de Prony's comments and other primary
documents indicate otherwise.

\subsection{\emph{Limites}}
\label{s11b}

Cauchy defined limits as follows in his \emph{Cours d'Analyse} (CdA):
\begin{quote}
On nomme quantit\'e \emph{variable} celle que l'on consid\`ere comme
devant recevoir successivement plusieurs valeurs diff\'erentes les
unes des autres.  \ldots{} Lorsque les valeurs successivement
attribu\'ees \`a une m\^eme variable s'app\-rochent ind\'efiniment
d'une valeur fixe, de mani\`ere \`a finir par en diff\'erer aussi peu
que l'on voudra, cette derni\`ere est appel\'ee la \emph{limite} de
toutes les autres.%
\footnote{Translation from \cite[p.\;6]{BS}: ``We call a quantity
\emph{variable} if it can be considered as able to take on
successively many different values.\;\ldots\;When the values
successively attributed to a particular variable indefinitely approach
a fixed value in such a way as to end up by differing from it by as
little as we wish, this fixed value is called the \emph{limit} of all
the other values.''}
(Cauchy \cite{Ca21}, 1821, p.\;4; emphasis in the original)
\end{quote}
Here Cauchy defines limits in terms of a primitive notion of a
\emph{variable quantity}.  As Robinson pointed out, Cauchy
``assign[ed] a central role to the notion of a variable which tends to
a limit, in particular to the limit zero'' (Robinson \cite{Ro66},
1966, p.\;276).

Elsewhere in CdA, Cauchy used what appears to be a somewhat different
notion of limit, as for example when the value of the derivative is
extracted from the ratio of infinitesimals~$\Delta y$ and~$\Delta x$
(see Section~\ref{s21}).  Two distinct approaches used by Cauchy are
analyzed in Section~\ref{s12}.

\subsection{A-track and B-track for the development of analysis}
\label{s12}

The article Katz--Sherry \cite{13f} introduced a distinction between
two types of procedures in the writing of the pioneers of
infinitesimal calculus:

\begin{enumerate}
\item[(A)] procedures in pioneering work in analysis that can be based
on an Archime\-dean continuum (or the A-track approach),
cf.\;\cite{Ar97}; and
\item[(B)] procedures that can be based on a Bernoullian (i.e.,
infinitesimal-enriched) continuum (the B-track approach), as they
appear in Leibniz, Bernoulli, Euler, and others.
\end{enumerate}

This is not an exhaustive distinction, but one that helps broaden the
lens of a historiography often handicapped by self-imposed limitations
of a Weierstrassian type; see Section~\ref{s12d}.

Here we use the term \emph{procedure} in a broad sense that
encompasses algorithms but is not limited to them.  For instance,
Euler's proof of the infinite product formula for the sine function is
a rather coherent procedure though it can hardly be described as an
algorithm; see \cite{17b} for an analysis of Euler's proof.

Like Leibniz, Cauchy used both A-track and B-track techniques in his
work.  The sample discussed in Section~\ref{s32} below illustrates his
A-track work.  Elsewhere, as we document in this article and in
earlier work (see e.g., \cite{12b}), Cauchy used B-track techniques,
as well.

\subsection{What is Cauchy's \emph{limite}?}
\label{s12c}

Scholars who stress Cauchy's use of the limit concept rely on a
traditional but flawed dichotomy of infinitesimals \emph{vs} limits.
The dichotomy is flawed because limits are present whether one works
with an Archimedean or Bernoullian continuum (see Section~\ref{s12}).
In fact, the definition of derivative found in Cauchy (see
Section~\ref{s21}) suggests that he works with the B-track version of
limits which is referred to as the \emph{standard part function} in
modern infinitesimal analysis; see Section\;\ref{s71},
formula\;\eqref{e83}.  Thus the real issue is whether Cauchy's
continuum was Archimedean or Bernoullian, and the genuine dichotomy is
between A-track~$\epsilon,\delta$ techniques and B-track infinitesimal
techniques.

\subsection{Butterfly model}
\label{s12d}

The articles (Bair et al.\;\cite{17a}), (Bair et al.\;\cite{17b}), and
(Fletcher et al.\;\cite{17f}) argued that some historians of
mathematics operate within a conceptual scheme described in (Hacking
\cite{Ha14}, 2014) as a \emph{butterfly model} of development.

Inspired in part by (Mancosu \cite{Ma09}, 2009), Ian Hacking proposes
a distinction between the \emph{butterfly model} and the \emph{Latin
model}, namely the contrast between a model of a deterministic
(genetically determined) biological development of animals like
butterflies (the egg--larva--cocoon--butterfly development), as
opposed to a model of a contingent historical evolution of languages
like Latin.

Historians working within the butterfly paradigm often assume that the
evolution of mathematical rigor has a natural direction, leading
\emph{forward} to the Archimedean framework as developed by
Weierstrass and others (what Boyer referred to as ``the great
triumvirate'' \cite[p.\;298]{Bo49}).  Such historians also tend to
interpret the qualifier \emph{rigorous} as necessarily implying
\emph{Archimedean}, as we illustrate in Section~\ref{s11c}.

\subsection{Siegmund-Schultze on \emph{Cours d'Analyse}}
\label{s11c}

As an illustration of butterfly model thinking by modern historians,
we turn to a review by historian Siegmund-Schultze of an English
edition of CdA (Bradley--Sandifer \cite{BS}, 2009).  The review
illustrates the poignancy of Grattan-Guinness' comment quoted in our
epigraph.  The comment appears in (Grattan-Guinness \cite{Gr70}, 1970,
p.\;379) in the context of a discussion of CdA.

Siegmund-Schultze's Zentralblatt (Zlb) review (\cite{Si09}, 2009) of
the English edition of CdA contains two items of interest:
\begin{enumerate}
[label={(SS\theenumi)}]
\item
\label{c1}
 Siegmund-Schultze quotes part of Cauchy's definition of continuity
via infinitesimals, and asserts that Cauchy's use of infinitesimals
was a step backward: ``There has been \ldots{}\;an intense historical
discussion in the last four decades or so how to interpret certain
apparent remnants of the past or -- as compared to J. L. Lagrange's
(1736--1813) rigorous `Algebraic Analysis' -- even \emph{steps
backwards} in Cauchy's book, particularly his use of
infinitesimals\ldots{}'' (\cite{Si09}; emphasis added).
\item
\label{c2}
Siegmund-Schultze quotes Cauchy's comments (in translation) on rigor
in geometry, and surmises that the framework for CdA was
Archi\-medean, similarly to Euclid's geometry: ``a non-Archi\-median
interpretation of the continuum would clash with the Euclidean theory,
which was still the basis of Cauchy's book. Indeed, Cauchy writes in
the `introduction' to the Cours d'Analyse: `As for methods, I have
sought to give them all the \emph{rigor that one demands in geometry},
\ldots'\,'' (ibid.; emphasis added).
\end{enumerate}
Siegmund-Schultze's Zbl review goes on to continue the quotation from
Cauchy:
\begin{quote}
``\ldots\;in such a way as never to revert to reasoning drawn from the
\emph{generality of algebra}.  Reasoning of this kind, although
commonly admitted, particularly in the passage from convergent to
divergent series and from real quantities to imaginary expressions,
can, it seems to me, only occasionally be considered as inductions
suitable for presenting the truth, since they accord so little with
the precision so esteemed in the mathematical sciences.''  (Cauchy as
quoted in \cite{Si09}; emphasis added).
\end{quote}
Cauchy's objections here have to do with the cavalier use of divergent
series, based on a heuristic principle Cauchy called the
\emph{generality of algebra}, by his illustrious predecessors Euler
and Lagrange, rather than with the issue of using or not using
infinitesimals, contrary to Siegmund-Schultze's claim.  We will
evaluate Siegmund-Schultze's claims further in Section~\ref{s14}.

\subsection{Analysis of a review}
\label{s14}

The Zbl review quoted in Section~\ref{s11c} tends to confirm the
diagnosis following Hacking.  Namely, the comment on infinitesimals
quoted in~\ref{c1} leading specifically \emph{backward} will surely be
read by the Zbl audience as indicative of an assumption of an organic
(butterfly model) \emph{forward} direction (culminating in the
\emph{great triumvirate}).  

Similarly, the comment quoted in~\ref{c2} appears to take it for
granted that Euclid's framework, being rigorous, was necessarily
Archimedean.  Yet the facts are as follows:
\begin{enumerate}
\item[(i)] Books I through IV of \emph{The Elements} are developed
without the Archimedean axiom;
\item[(ii)] developments around 1900 showed conclusively that the
completeness property of~$\R$ is irrelevant to the development of
Euclidean geometry, and in fact the latter can be developed in the
context of non-Archimedean fields.  
\end{enumerate}
Indeed, Hilbert proved that these parts of Euclidean geometry can be
developed in a non-Archimedean plane (modulo some specific assumptions
such as circle--circle intersection and postulation of the congruence
theorems); see further in \cite[Section\;5]{13a}.  

While Euclid relied on the Archimedean axiom to develop his
\emph{theory of proportion}, Hilbert obtained all the results of
Euclidean geometry including the theory of proportion and geometric
similarity without such a reliance; see Hartshorne (\cite{Ha00}, 2000,
Sections 12.3--12.5 and 20--23) or Baldwin (\cite{Ba17}, 2017).

Furthermore, starting with Descartes' \emph{Geometry}, mathematicians
implicitly relied on ordered field properties rather than the ancient
theory of proportion.

Moreover, it is difficult to understand how Siegmund-Schutze would
reconcile his two claims.  If Cauchy used Euclidean Archimedean
mathematics exclusively, as implied by \ref{c2}, then what exactly
were the entities that constituted a step backward, as claimed in
\ref{c1}?  Siegmund-Schultze's counterfactual claims are indicative of
butterfly-model thinking as outlined in Section~\ref{s12d}.

Like the Zbl review by Siegmund-Schultze, the Cauchy scholarship of
Gilain tends to be colored by teleological assumptions of the sort
detailed above, as we argue in Sections~\ref{s2} and \ref{s5}.  

A number of historians and mathematicians have sought to challenge the
received views on Cauchy's infinitesimals, as we detail in
Sections~\ref{s1} through \ref{laug}.

\subsection{Robinson on received views}
\label{s1}

Abraham Robinson noted that the received view of the development of
the calculus
\begin{quote}
[would] lead us to expect that, following the rejection of Leibniz'
theory by Lagrange and D'Alembert, infinitely small and infinitely
large quantities would have no place among the ideas of Cauchy, who is
generally regarded as the founder of the modern approach, or that they
might, at most, arise as figures of speech, as in `$x$\;tends to
infinity'.  However, this expectation is mistaken.
\cite[p.\;269]{Ro66}.
\end{quote}
Robinson described Cauchy's approach as follows:
\begin{quote}
Cauchy regarded his theory of infinitely small quantities as a
satisfactory foundation for the theory of limits and (d'Alembert's
suggestion notwithstanding) he \emph{did not introduce the latter in
order to replace the former}.  His proof procedures thus involved
\emph{both} infinitely small (and infinitely large) quantities
\emph{and} limits.  \cite[p.\;271]{Ro66} (emphasis added)
\end{quote}
Note Robinson's focus on Cauchy's \emph{procedures} (for a discussion
of the procedure/ontology dichotomy, see B\l aszczyk et
al.\;\cite{17d}).  After quoting Cauchy's definition of derivative,
Robinson notes:
\begin{quote}
Later generations have overlooked the fact that in this definition
$\Delta x$ and~$\Delta y$ were explicitly supposed to be infinitely
small.  Indeed according to our present standard ideas, we take
$f'(x)$ to be the limit [of]~$\Delta y/\Delta x$ as~$\Delta x$ tends
to zero, whenever that limit exists, without any mention of infinitely
small quantities.  Thus, as soon as we consider limits, the assumption
that~$\Delta x$ and~$\Delta y$ are infinitesimal is completely
redundant.  It is therefore the more interesting that the assumption
is there, and, indeed, appears again and again also in Cauchy's later
expositions of the same topic (Cauchy [1829, 1844]).
\cite[p.\;274]{Ro66}
\end{quote}
Robinson's conclusion is as follows:
\begin{quote}
We are forced to conclude that Cauchy's mental picture of the
situation was significantly different from the picture adopted today,
in the Weierstrass tradition.  (ibid.)
\end{quote}
It is such received views in what Robinson refers to as the
\emph{Weierstrass tradition} that we wish to reconsider here.

\subsection{Grattan-Guinness on Cauchy's infinitesimals}
\label{s11}

Robinson's 1966 comments on the Weierstrassian tradition cited in
Section~\ref{s1} were echoed by historians Ivor Grattan-Guinness and
Detlef Laugwitz.  Thus, fourteen years later, Grattan-Guinness wrote:
\begin{quote}
[Cauchy's definition of infinitesimal] is \emph{in contrast to} the
view adopted from the Weierstrassians onwards (and occasionally
earlier), where an infinitesimal is a variable with \emph{limit}
zero\ldots{} (Grattan-Guinness \cite{Gr80}, 1980, p.\;110; emphasis
added)
\end{quote}
Concerning the term \emph{limit}, it is necessary to disassociate the
following two issues:
\begin{enumerate}
[label={(Ca\theenumi)}]
\item
\label{i1}
the issue of whether or not limits were at the base of Cauchy's
approach;
\item
\label{i2}
the issue of Cauchy's systematic use of infinitesimals as numbers in
his textbooks and research articles.
\end{enumerate}

\subsection{Laugwitz on Cauchy's infinitesimals}
\label{laug}

As far as item~\ref{i2} is concerned, Laugwitz acknowledged that
Cauchy started using infinitesimals systematically in the 1820s
(whereas his attitude toward them during the preceding decade was more
ambiguous and limits may have played a larger role):
\begin{quote}
\ldots \emph{after 1820, Cauchy developed his analysis by utilizing
infinitesimals in a deliberate and consequent manner.}  (Laugwitz
\cite{La89}, 1989, p.\;196; emphasis in the original)
\end{quote}
Laugwitz' position is consistent with Gilain's observation that
infinitesimals first appeared in Cauchy's course summary during the
academic year 1820--1821:
\begin{quote}
Ann\'ee 1820--1821 \ldots\;Notons aussi l'apparition, pour la
premi\`ere fois dans les \emph{Mati\`eres des le\c{c}ons}, des notions
de quantit\'es infiniment petites et infiniment grandes
(le\c{c}on\;3).%
\footnote{Translation: ``Year 1820--1821 \ldots\;We also note the
appearance, for the first time in the \emph{Lesson summaries}, of the
notions of infinitely small and infinitely large quantities (lesson
3).''}
(Gilain \cite{Gi89}, \S 52, 1989)
\end{quote}
In 1997, Laugwitz elaborated on the subject (of Cauchy's endorsement
of infinitesimals circa 1820) in the following terms:
\begin{quote}
Cauchy avoided the use of the infinitely small. This provoked growing
criticism on the part of his colleagues, including the physicist
Petit, who emphasized the didactical and practical advantages of the
use of infinitely small magnitudes. In 1819 and in 1820, the Conseil
d'Instruction at the Ecole exerted strong pressure on Cauchy, but this
alone would not have made this rather stubborn man change his mind.
\emph{Around 1820, he must have realized that infinitesimal
considerations were a powerful research method} at a time when he was
in a state of constant rivalry, especially with Poisson.  (Laugwitz
\cite{La97}, 1997, p.\;657; emphasis added)
\end{quote}
In the textbook \emph{Cours d'Analyse} \cite{Ca21}, \emph{limite} is
not the only central foundational concept for Cauchy, as we argue in
Section~\ref{s2}.

We challenge a common misconception according to which Cauchy was
forced to teach infinitesimals at the \emph{Ecole Polytechnique}
allegedly against his will.  We show that the debate there concerned
mainly the issue of \emph{rigor}, a separate one from
\emph{infinitesimals}; see Section~\ref{s5}.

\section{Cauchy's \emph{limite} and \emph{infiniment petit}}
\label{s2}

In this section we will analyze the meaning of Cauchy's terms
\emph{limite} and \emph{infiniment petit}.

\subsection{Differentials and infinitesimals}
\label{s21}

In his work, Cauchy carefully distinguishes between
differentials~$ds,dt$ which to Cauchy are noninfinitesimal variables,
on the one hand, and infinitesimal increments~$\Delta s, \Delta t$, on
the other:
\begin{quote}
\ldots soit~$s$ une variable distincte de la variable primitive~$t$.
En vertu des d\'efinitions adopt\'ees, le rapport entre les
diff\'erentielles~$ds, dt$, sera la limite du rapport entre les
\emph{accroissements infiniment petits}~$\Delta s, \Delta t$.%
\footnote{Translation: ``Let $s$ be a variable different from the
primitive variable~$t$.  By virtue of the definitions given, the ratio
of the differentials $ds, dt$ will be the limit of the ratio of the
infinitely small increments $\Delta s, \Delta t$.''}
(Cauchy \cite{Ca44}, 1844, p.\;11; emphasis added)
\end{quote}
Cauchy goes on to express such a relation by means of a formula in
terms of the infinitesimals~$\Delta s$ and~$\Delta t$:
\begin{quote}
On aura donc
\begin{equation}
\label{e21}
\frac{ds}{dt} = \; \text{lim.}\, \frac{\Delta s}{\Delta t}
\end{equation}
(ibid., equation (1); the period after lim in ``lim.'' in the
original; equation number~\eqref{e21} added)
\end{quote}
Cauchy's procedure involving the passage from the ratio of
infinitesimals like~$\frac{\Delta s}{\Delta t}$ to the value of the
derivative $\frac{ds}{dt}$ as in equation~\eqref{e21} has a close
parallel in Robinson's infinitesimal analysis, where it is carried out
by the standard part function; see equations~\eqref{e61}
and~\eqref{e62} in Section~\ref{s71}.

Paraphrasing this definition in Archimedean terms would necessarily
involve elements that are not explicit in Cauchy's definition.  Thus
Cauchy's ``lim.''\;finds a closer proxy in the notion of standard
part, as in formula~\eqref{e83}, than in any notion of limit in the
context of an Archimedean continuum; see also~Bascelli et
al.\;(\cite{14a}, 2014).

\subsection{Definite integrals and infinitesimals}
\label{s22c}

Similar remarks apply to Cauchy's 1823 definition of the definite
integral which exploits a partition of the domain of integration into
infinitesimal subintervals.  Here Cauchy writes: ``D'apr\`es ce qui a
\'et\'e dit dans la derni\`ere le\c{c}on, si l'on divise~$X-x_0$ en
\'el\'emens%
\footnote{We preserved the original spelling.}
infiniment petits~$x_1-x_0, x_2-x_1\,\ldots\,X-x_{n-1}$, la somme
\[
(1) \quad S =
    (x_1-x_0)f(x_0)+(x_1-x_2)f(x_1)+\ldots+(X-x_{n-1})f(x_{n-1})
\]
convergera vers une \emph{limite} repr\'esent\'ee par l'int\'egrale
d\'efinie
\[
(2) \hskip1.9in \int_{x_0}^X f(x)dx. \hskip1.9in
\]
Des principes sur 
%
%spelling of lesquels corrected by refereeA for antiquitates
%
lesquels nous avons fond\'e cette proposition il r\'esulte, etc.''
(Cauchy \cite{Ca23}, 1823, Le\c{c}on~22, p.\;85; emphasis added).

Note that there is a misprint in Cauchy's formula~(1): the difference
$(x_1-x_2)$ should be~$(x_2-x_1)$.  In this passage, Cauchy refers to
the successive differences~$x_1-x_0$,~$x_2-x_1$,~$X-x_{n-1}$ as
\emph{infinitely small elements}.

Analogous partitions into infinitesimal subintervals are exploited in
Keisler's textbook \cite{Ke86} (and throughout the literature on
infinitesimal analysis; see e.g., \cite[p.\;153]{Go98}).  Cauchy's use
of \emph{limite} in the passage above is another instance of limit in
the context of a Bernoullian continuum, which parallels the use of the
standard part function (see Section~\ref{s71}) enabling the transition
from a sum of type~(1) above to the definite integral~(2), similar to
the definition of the derivative analyzed in Section~\ref{s21}.

\subsection{\emph{Un infiniment petit} in Cauchy}
\label{s23}

What is the precise meaning of Cauchy's \emph{infiniment petit}
(infinitely small)?  All of Cauchy's textbooks on analysis contain
essentially the same definition up to slight changes in word order:
\begin{quote}
Lorsque les valeurs num\'eriques successives d'une m\^eme variable
d\'ecroissent ind\'efiniment, de mani\`ere \`a s'abaisser au-dessous
de tout nombre donn\'e, cette variable devient ce qu'on nomme un
\emph{infiniment petit} ou une quantit\'e \emph{infiniment petite}.
Une variable de cette esp\`ece a z\'ero pour limite.%
\footnote{Translation: ``When the successive numerical values of such
a variable decrease indefinitely, in such a way as to fall below any
given number, this variable becomes what we call \emph{infinitesimal},
or an \emph{infinitely small quantity}.  A variable of this kind has
zero as its limit'' \cite[p.\;7]{BS}.}
\cite[p.\;4]{Ca21} (emphasis in the original)
\end{quote}
An examination of the books \cite{Ca21}, \cite{Ca23} reveals that
Cauchy typically \emph{did not} define his infinitely small literally
as a variable whose limit is zero.  Namely, he rarely wrote ``an
infinitely small \emph{is} a variable, etc.''  but said, rather, that
a variable \emph{becomes} (\emph{devient}) an infinitely small.

Thus, the passage cited above is the first definition of the
infinitely small in \emph{Cours d'Analyse}.  The next occurrence is on
page\;26 there, again using \emph{devient}, and emphasizing
\emph{infiniment petite} by means of italics.  On page\;27 Cauchy
summarizes the definition as follows: ``Soit~$\alpha$ une quantit\'e
infiniment petite, c'est-\`a-dire, une variable dont la valeur
num\'erique d\'ecroisse ind\'efiniment.''  This is a summary of the
definition already given twice, the expression ``infiniment petite''
is not italicized, and ``is'' is used in place of ``becomes'' as
shorthand for the more detailed and precise definitions appearing
earlier in Cauchy's textbook.  An identical definition with
\emph{devient} appears in his 1823 textbook \cite[p.\;4]{Ca23}.

Cauchy's term \emph{becomes} implies a change of nature or
\emph{type}.%
\footnote{To illustrate such a change in modern terms, note that in
the context of the traditional construction of the real numbers in
terms of Cauchy sequences~$u=(u_n)\in\Q^{\N}$ of rational numbers, one
never says that a real number \emph{is} a sequence, but rather that a
sequence \emph{represents} or \emph{generates} the real number, or to
use Cauchy's terminology, \emph{becomes} a real number.  A related
construction of hyperreal numbers out of sequences of real numbers,
where a sequence tending to zero generates an infinitesimal, is
summarized in Section~\ref{s71}.}
Namely, a variable is not quite an infinitesimal yet, but only serves
to \emph{generate} or \emph{represent} one, as emphasized by Laugwitz:
\begin{quote}
Cauchy never says what his infinitesimals \emph{are}; we are told only
how infinitesimals can be \emph{represented}.  (Laugwitz \cite{La87},
1987, p.\;271)
\end{quote}
See also Sad et al.\;\cite{Sa01}.  This indicates that Cauchy
considered an infinitesimal as a separate type of mathematical entity,
distinct from variable or sequence.

\subsection{Variable quantities, infinitesimals, and limits}
\label{s24}

To comment more fully on Cauchy's passage cited in Section~\ref{s23},
note that there are three players here:
\begin{enumerate}
\item[(A)] variable quantity;
\item[(B)] infinitesimal;
\item[(C)] limit zero.
\end{enumerate}
We observe that the notion of variable quantity is the primitive
notion in terms of which both infinitesimals and limits are defined
(see Section~\ref{s11b} for Cauchy's definition of limit in terms of
variable quantity).  This order of priorities is confirmed by the
title of Cauchy's very first lesson in his 1823 book:
\begin{quote}
1.$^{re}$ Le\c con.  Des variables, de leurs limites, et des
quantit\'es infiniment petites \cite[p.\;ix]{Ca23}
\end{quote}
Thus, Cauchy is proposing a definition and an observation:
\begin{enumerate}
[label={(Co\theenumi)}]
\item
\label{i1c}
a variable quantity that diminishes indefinitely becomes an
infinitesimal; and
\item
\label{i2c}
such a variable quantity has zero as limit.
\end{enumerate}
Here item~\ref{i2c} is merely a restatement of the property of
diminishing indefinitely in terms of the language of limits.  As noted
in Section~\ref{s1b}, Robinson pointed out that Cauchy assigned a
central role to the notion of a variable which tends to a limit.
Cauchy's notion of limit here is close to the notion of limit of his
predecessor Lacroix (see Section~\ref{s22}).

\subsection{Assigning a sign to an infinitesimal}

Cauchy often uses the notation~$\alpha$ for a generic infinitesimal,
in both his 1821 and 1823 textbooks.  In his 1823 textbook Cauchy
assumes that~$\alpha$ is either positive or negative:
\begin{quote}
Cherchons maintenant la limite vers laquelle converge l'expression
$(1+ \alpha)^{\frac{1}{\alpha}}$, tandis que~$\alpha$ s'approche
ind\'efini\-ment de z\'ero.  Si l'on suppose d'abord la quantit\'e
$\alpha$ positive et de la forme~$\frac{1}{m}$,~$m$ d\'esignant un
nombre entier variable et susceptible d'un accroissement ind\'efini,
on aura~$(1+\alpha)^{\frac{1}{\alpha}}=\left(1+\tfrac{1}{m}\right)^m$
\ldots{} Supposons enfin que~$\alpha$ devienne une quantit\'e
n\'egative.  Si l'on fait dans cette hypoth\`ese
$1+\alpha=\frac{1}{1+\beta}$, \,$\beta$ sera une quantit\'e positive,
qui convergera elle-m\^eme vers z\'ero, etc.  \cite[pp.\;2--4]{Ca23}
\end{quote}
It is well known that variable quantities or sequences that generate
Cauchyan infinitesimals are not necessarily monotone.  Indeed, Cauchy
himself gives a non-monotone example at the beginning of CdA:
\begin{quote}
$\frac14, \frac13, \frac16, \frac15, \frac18, \frac17,$ \&c.\,\ldots{}
\cite[p.\;27]{Ca21}
\end{quote}
This poses a problem since it is not obvious how to assign a sign plus
or minus to an arbitrary null sequence (i.e., a sequence tending to
zero).

When Cauchy actually uses infinitesimals in proofs and applications,
he assumes that they can be manipulated freely in arithmetic
operations and other calculations.  While formal order theory is a few
decades away and is not to be found as such in Cauchy, he does appear
to assume that a definite sign can be attached to an infinitesimal.
Besides assuming that they have a well-defined sign, Cauchy also
routinely applies arithmetic operations to infinitesimals.

This creates a difficulty to those who consider that Cauchy merely
used the term ``infinitely small'' as shorthand for a sequence with
limit~$0$, since it is unclear how to assign a sign to an arbitrary
null sequence, whereas Cauchy does appear to assign a sign to his
infinitesimals.

Which process exactly did Cauchy envision when he spoke of a sequence
\emph{becoming} an infinitesimal?  Cauchy does not explain.  However,
Cauchy's assumption that each infinitesimal has a sign suggests that a
sequence is not identical to the infinitesimal it generates.

Even monotone sequences are not closed under arithmetic operations.
Namely, such operations necessarily lead to non-monotone ones,
including ones that change sign.

Cauchy routinely assumes in his work, particularly on integrals, that
one can freely add infinitesimals and obtain other infinitesimals,
i.e., that the numbers involved are closed under arithmetic
operations.

Such an assumption is valid in modern theories of ordered fields
properly extending~$\R$, but if one is working with sequences, such an
assumption leads to a dilemma:
\begin{enumerate}
\item
either one only works with monotone ones, in which case one gets into
a problem of closedness under natural arithmetic operations;
\item
or one works with arbitrary sequences, in which case the
assumption that a sequence can be declared to be either positive or
negative becomes problematic.
\end{enumerate}
Cauchy was probably not aware of the difficulty that that one can't
\emph{both} assign a specific sign to~$\alpha$, and also have the
freedom of applying arithmetic operations to infinitesimals.  The
point however is that the way he uses infinitesimals indicates that
both conditions are assumed, even though from the modern standpoint
the justification provided is insufficient.  In other words, Cauchy's
\emph{procedures} are those of an infinitesimal-enriched framework,
though the \emph{ontology} of such a system is not provided.

Cauchy most likely was not aware of the problem, for otherwise he may
have sought to address it in one way or another.  He did have some
interest in asymptotic behavior of sequences.  Thus, in some of his
texts from the late 1820s he tried to develop a theory of the order of
growth at infinity of functions.  Such investigations were eventually
picked up by du Bois-Reymond, Borel, and Hardy; see Borovik--Katz
(\cite{12b}, 2012) for details.

\subsection
{Gilain on omnipresence of limits}
\label{s26}

Gilain refers to Cauchy's course in 1817 as\;a
\begin{quote}
cours tr\`es important historiquement, o\`u les bases de la nouvelle
analyse, notamment celle de l'\emph{Analyse alg\'ebrique} de 1821,
sont pos\'ees\ldots\;\cite[\S 30]{Gi89}
\end{quote}
He goes on to note ``l'omnipr\'esence du concept de limite'' (ibid.).
How are we to evaluate Gilain's claim as to the ``omnipresence'' of
the concept of limit?

With regard to Cauchy's pre-1820 courses such as the one in 1817
mentioned by Gilain, there appears to be a consensus among scholars
already noted in Section~\ref{s11} concerning the absence of
infinitesimals.  As far as Cauchy's 1821 book is concerned, the
presence (perhaps even ``omnipresence'' as per Gilain) of limits in
the definition of infinitesimals goes hand-in-hand with the fact that
Cauchy defined both limits and infinitesimals in terms of the
primitive notion of a \emph{variable quantity} (see beginning of
Section~\ref{s1b} as well as Section~\ref{s24}).  It is therefore
difficult to agree with Gilain when he claims to know the following:
\begin{quote}
On sait que Cauchy d\'efinissait le concept d'infiniment petit \`a
l'aide du concept de limite, qui avait le premier r\^ole (voir Analyse
alg\'ebrique, p.\;19; \ldots) \cite[note~67]{Gi89}
\end{quote}
Here Gilain claims that it is the concept of \emph{limite} that played
a primary role in the definition of infinitesimal, with reference to
page 19 in the 1897 \emph{Ouevres Compl\`etes} edition of CdA
\cite{Ca21}.  The corresponding page in the 1821 edition is page 4.
We quoted Cauchy's definition in Section~\ref{s23} and analyzed it in
Section~\ref{s24}.  An attentive analysis of the definition indicates
that it is more accurate to say that it is the concept of variable
quantity (rather than \emph{limite}) that ``avait le premier r\^ole.''

Cauchy exploited the notion of limit in \cite[Chapter~2, \S 3]{Ca21}
in the proofs of Theorem~1 and Theorem~2.  Theorem~1 compares the
convergence of the difference~$f(x+1)-f(x)$ and that of the
ratio~$\frac{f(x)}{x}$.  Theorem~2 compares the convergence
of~$\frac{f(x+1)}{f(x)}$ and~$[f(x)]^{\frac{1}{x}}$.  These proofs can
be viewed as prototypes of~$\epsilon,\delta$ arguments.  On the other
hand, neither of the two proofs mentions infinitesimals.  Therefore
neither can support Gilain's claim to the effect that Cauchy allegedly
used limits as a basis for defining infinitesimals.  The proof of
Theorem\;1 is analyzed in more detail in Section~\ref{s32}.

Cauchy's procedures exploiting infinitesimals have stood the test of
time and proved their applicability in diverse areas of mathematics,
physics, and engineering.

Gilain and some other historians assume that the appropriate modern
proxy for Cauchy's \emph{limite} necessarily operates in the context
of an Archimedean continuum (see Section~\ref{s24}).  Yet the
vitality and robustness of Cauchy's infinitesimal procedures is
obvious given the existence of proxies in modern theories of
infinitesimals.  What we argue is that modern infinitesimal proxies
for Cauchy's procedures are more faithful to the original than
Archimedean proxies that typically involve anachronistic paraphrases
of Cauchy's briefer definitions and arguments.

This article does not address the historical \emph{ontology} of
infinitesimals (a subject that may require separate study) but rather
the \emph{procedures} of infinitesimal calculus and analysis as found
in Cauchy's oeuvre (see~\cite{17d} for further details on the
procedure/ontology dichotomy).

\subsection{\emph{Limite} and infinity}
\label{s12b}

As we noted in Section~\ref{s12c}, the use of the term \emph{limite}
by Cauchy could be misleading to a modern reader.  Consider for
example its use in the passage cited in Section\;\ref{s23}.  The fact
that Cauchy is not referring here to a modern notion of limit is
evident from his very next sentence:
\begin{quote}
Lorsque les valeurs num\'eriques successives d'une m\^eme variable
croissent de plus en plus, de mani\`ere \`a s'\'elever au-dessus de
tout nombre donn\'e, on dit que cette variable a pour limite l'infini
positif indiqu\'e par le signe~$\infty$ s'il s'agit d'une variable
positive\ldots%
\footnote{Translation: ``When the successive numerical values [i.e.,
absolute values] of the same variable grow larger and larger so as to
rise above each given number, one says that this variable has limit
positive infinity denoted by the symbol~$\infty$ when the variable is
positive.''}
\cite[p.\;4]{Ca23}
\end{quote}
In today's calculus courses, it is customary to give an
$(\epsilon,\delta)$ or~$(\epsilon,N)$ definition of limit of, say, a
sequence, and then introduce infinite `limits' in a broader sense when
the sequence diverges to infinity.  But Cauchy does not make a
distinction between convergent limits and divergent infinite limits.

Scholars ranging from Sinaceur (\cite{Si73}, 1973) to Nakane
(\cite{Na14}, 2014) have pointed out that Cauchy's notion of limit is
distinct from the Weierstrassian \emph{Epsilontik} one (this is
particularly clear from Cauchy's definition of the derivative analyzed
in Section~\ref{s21}); nor did Cauchy ever give an~$\epsilon,\delta$
\emph{definition} of limit, though prototypes of~$\epsilon,\delta$
\emph{arguments} do occasionally appear in Cauchy; see
Section~\ref{s12}.

\section{Minutes of meetings, Poisson, and de Prony}
\label{s5}

Here we develop an analysis of the third of the misconceptions
diagnozed in Borovik--Katz (\cite{12b}, 2012, Section 2.5), namely the
idea that Cauchy was forced to teach infinitesimals at the \emph{Ecole
Polytechnique} allegedly against his will.  We show that the debate
there concerned mainly the issue of \emph{rigor}, a separate one from
\emph{infinitesimals}.

Minutes of meetings at the \emph{Ecole} are a valuable source of
information concerning the scientific and pedagogical interactions
there in the 1820s.

\subsection{Cauchy pressured by Poisson and de Prony}

Gilain provides detailed evidence of the pressure exerted by Sim\'eon
Denis Poisson, Gaspard de Prony, and others on Cauchy to simplify his
analysis course.  Thus, in 1822
\begin{quote}
Poisson et de Prony\ldots{} insistent [sur la] n\'ecessit\'e\ldots{}
de simplifier l'enseignement de l'analyse, en multipliant les exemples
num\'eriques et en r\'eduisant beaucoup la partie analyse alg\'ebrique
plac\'ee au d\'ebut du cours.  \cite[\S 61]{Gi89}
\end{quote}
Similarly, in 1823, Cauchy's course was criticized for being too
complicated:
\begin{quote}
des voix se sont \'elev\'ees pour trouver trop compliqu\'ees les
feuilles de cours en question et il \'etait d\'ecid\'e de proposer au
Ministre la nomination d'une commission qui serait charg\'ee chaque
ann\'ee de l'examen des feuilles d'analyse et des modifications
\'eventuelles \`a y apporter.  \cite[\S 72]{Gi89}
\end{quote}
The critics naturally include Poisson and de Prony:
\begin{quote}
Cette commission, effectivement mise en place, comprendra, outre
Laplace, pr\'esident, les examinateurs de math\'ematiques (Poisson et
de Prony),\ldots{} (ibid.)
\end{quote}
The complaints continue in 1825 as Fran\c cois Arago declares that
\begin{quote}
ce qu'il y a de plus utile \`a faire pour le cours d'analyse, c'est de
le simplifier.  \cite[\S 84]{Gi89}
\end{quote}
At this stage Cauchy finally caves in and declares (in third person):
\begin{quote}
il ne s'attachera plus \`a donner, comme il a fait jusqu'\`a
pr\'esent, des d\'emonstrations parfaitement \emph{rigoureuses}.
\cite[\S 86]{Gi89} (emphasis added)
\end{quote}
Note however that in these discussions, the issue is mainly that of
\emph{rigor} (i.e., too many proofs) rather than choice of a
particular approach to the foundations of analysis.  While Cauchy's
commitment to simplify the course may have entailed skipping the
proofs in the style of the \emph{Epsilontik} of Theorems~1 and 2 in
\cite[Chapter 2, \S 3]{Ca21} (see end of Section~\ref{s24}), it may
have also entailed skipping the proofs of as many as \emph{eight}
theorems concerning the properties of \emph{infinitesimals} of various
orders in \cite[Chapter 2, \S 1]{Ca21}, analyzed in
\cite[Section~2.3]{12b}.

\subsection{Reports by de Prony}
\label{s36}

Gilain notes that starting in 1826, there is a new source of
information concerning Cauchy's course, namely the reports by de
Prony:
\begin{quote}
de Prony reproche de fa\c con g\'en\'erale \`a Cauchy de \emph{ne pas
utiliser suffisamment} les consid\'erations g\'eom\'etriques et les
infiniment petits, tant en analyse qu'en m\'ecanique.  \cite[\S
101]{Gi89} (emphasis added)
\end{quote}
Thus with regard to the post-1820 period, only starting in 1826 do we
have solid evidence that not merely excessive \emph{rigor} but also
insufficient use of \emph{infinitesimals} was being contested.  Even
here, the complaint is not an alleged \emph{absence} of
infinitesimals, but merely \emph{insufficient use} thereof.  We will
examine de Prony's views in Section~\ref{s37}.

\subsection{Course summaries}

According to course summaries reproduced in \cite{Gi89}, Cauchy taught
both continuous functions and infinitesimals (and presumably the
definition of continuity in terms of infinitesimals after 1820) in the
\emph{premi\`ere ann\'ee} during the academic years 1825--1826,
1826--1827, 1827--1828, and 1828--1829 (the summaries for the
\emph{premi\`ere ann\'ee} during the 1829--1830 academic year,
Cauchy's last at the \emph{Ecole Polytechnique}, are not provided).
All these summaries contain identical comments on continuity and
infinitesimals for those years:
\begin{quote}
Des fonctions en g\'en\'eral, et des fonctions \emph{continues} en
particulier. -- Repr\'esentation g\'eom\'etrique des fonctions
\emph{continues} d'une seule variable. -- Du rapport entre
l'accroisse\-ment d'une fonction et l'accroissement de la variable.
-- Valeur que prend ce rapport quand les accroissemens deviennent
\emph{infiniment petits}.  (Cauchy as quoted by Gilain; emphasis
added)
\end{quote}
In 1827 for the first time we find a claim of an actual \emph{absence}
of infinitesimals from Cauchy's teaching.  Thus, on 12 january 1827,
\begin{quote}
le cours de Cauchy a de nouveau \'et\'e mis en cause pour sa
difficult\'e, (le gouverneur affirmant que des \'el\`eves avaient
d\'eclar\'e qu'ils ne le comprenaient pas), et son \emph{non-usage} de
la m\'ethode des infiniment petits (voir document C12).%
\footnote{To comment on Gilain's ``document C12'' (denoted C$_{12}$ in
\cite{Gi89}), it is necessary to reproduce what the document actually
says: ``Un membre demande si le professeur expose la m\'ethode des
infiniment petits, ainsi que le voeu en a \'et\'e exprim\'e.''  What
was apparently Cauchy's response to this query is reproduced in the
next paragraph of document C12: ``On r\'epond que le commencement du
cours ne pourra \^etre fond\'e sur les notions infinit\'esimales que
l'ann\'ee prochaine, parce que le cours de cette ann\'ee \'etait
commenc\'e \`a l'\'epoque o\`u cette disposition a \'et\'e
arr\^et\'ee; que M. Cauchy s'occupe de la r\'edaction de ses feuilles,
en cons\'equence, et qu'il a promis de les communiquer bient\^ot \`a
la commission de l'enseignement math\'ematique.''  

Thus, the actual contents of document C12 indicate that Gilain's claim
of ``\emph{non-usage}'' is merely an extrapolation.}
\cite[\S 103]{Gi89} (emphasis added)
\end{quote}
Tellingly, this comment by Gilain is accompanied by a footnote\;111
where Gilain acknowledges that in the end Cauchy did use
infinitesimals that year in his treatment of the theory of contact of
curves; see Section~\ref{s34} for details.

\subsection{Cauchy taken to task}
\label{s34}

Gilain writes that during the 1826--1827 academic year, Cauchy was
taken to task in the \emph{Conseil de Perfectionnement} of the
\emph{\'Ecole Polytechnique} for allegedly not teaching infinitesimals
(see \cite[\S103]{Gi89}).  Gilain goes on to point out in his
footnote\;111 that Cauchy exploited infinitesimals anyway that year,
in developing the theory of contact of curves:
\begin{quote}
S'il ne fonde pas le calcul diff\'erentiel et int\'egral sur la
`m\'ethode' des infiniment petits, Cauchy n'en utilise pas moins de
fa\c con importante ces objets (consid\'er\'es comme des variables
dont la limite est z\'ero),%
\footnote{Gilain's parenthetical remark here is an editorial comment
for which he provides no evidence.  The remark reveals more about
Gilain's own default expectations (see Section~\ref{s1b}) than about
Cauchy's actual foundational stance.}
en liaison notamment avec l'exposition de la th\'eorie du contact des
courbes.  \cite[note\;111]{Gi89}
\end{quote}
It emerges that Cauchy did use infinitesimals that year in his
treatment of a more advanced topic (theory of contact).  Thus Cauchy's
actual scientific practice was not necessarily dependent on his
preliminary definitions.  There is conflicting evidence as to whether
Cauchy used infinitesimals (as developed in \cite{Ca21} and
\cite{Ca23}) in the introductory part of his course that year.  As we
mentioned in Section~\ref{s36}, the course summary for 1826--1827 does
include both continuity and infinitesimals.

\subsection{Critique by de Prony}
\label{s37}

Gilain describes de Prony's criticism of Cauchy as follows:

\begin{quote}
[De Prony] critique notamment l'emploi de la m\'ethode des limites par
Cauchy au lieu de celle des infiniment petits, faisant appel ici \`a
l'autorit\'e posthume de Laplace, d\'ec\'ed\'e depuis le 5 mars 1827
(voir document C14). \\ \cite[\S 105]{Gi89}
\end{quote}
Here Gilain is referring to the following comments by de Prony:
\begin{quote}
Les d\'emonstrations des formules generales%
\footnote{\label{f9}The spelling as found in (Gilain
\cite[Document\;C$_{14}$]{Gi89}) is \emph{g\'en\'erales} (i.e., the
modern French spelling).  Gilain similarly replaced \emph{encor} by
\emph{encore}, \emph{mel\'ees} by \emph{m\^el\'ees},
\emph{immediatement} by \emph{imm\'ediatement}, \emph{methode} by
\emph{m\'ethode}, \emph{abrege} by \emph{abr\`ege}, and
\emph{collegue} by \emph{coll\`egue}.}
du mouvement vari\'e se sont encor trouv\'ees mel\'ees de
considerations relatives aux \emph{limites}; \ldots\;(de Prony as
quoted in Grattan-Guinness \cite{GG}, 1990, p.\;1339; emphasis in the
original)
\end{quote}
Having specified the target of his criticism, namely Cauchy's concept
of \emph{limite}, de Prony continues:
\begin{quote}
\ldots{}\;il me semble qu'en employant, immediatement et
exclusivement, la methode des infiniment petits, on abrege et on
simplifie les raisonnements sans nuire \`a la clart\'e; rappellons
nous combien cette methode \'etait recommand\'ee par l'illustre
collegue [Laplace] que la mort nous a enlev\'e.  (ibid.)
\end{quote}
What is precisely the nature of de Prony's criticism of Cauchy's
approach to analysis?  Does his criticism focus on excessive rigor, or
on infinitesimals, as Gilain claims?  The answer depends crucially on
understanding de Prony's own approach, explored in Section~\ref{s311}.

\subsection{De Prony on small oscillations}
\label{s311}

In his work \emph{M\'eca\-nique philo\-so\-phique}, de Prony considers
infinitesimal oscillations of the pendulum (de Prony \cite{De}, 1799,
p.\;86, \S125).  He gives the familiar formula for the period or more
precisely halfperiod, namely
\[
\pi\sqrt{\tfrac{a}{g}}
\]
where~$a$ is the length of the cord, and~$g$ is acceleration under
gravity.  Limits are not mentioned.  In the table on the following
page 87, he states the property of \emph{isochronism}, meaning that
the halfperiod~$\pi\sqrt{\frac{a}{g}}$ is independent of the size of
the infinitesimal amplitude.  This however is not true literally but
only up to a passage to limits, or taking the standard part;%
\footnote{Even if literally infinitesimal amplitudes are admitted,
there is still a discrepancy disallowing one to claim that the
halfperiod is literally~$\pi\sqrt{\frac{a}{g}}$.  This difficulty can
be overcome in the context of modern infinitesimal analysis; see
Kanovei et al.\;(\cite{16c}, 2016).}
see Section~\ref{s71}.  Thus de Prony's own solution to the conceptual
difficulties involving limits/standard parts in this case is merely to
ignore the difficulties and suppress the limits.

In his article ``Suite des le\c cons d'analyse,'' de Prony lets~$n=Az$
(\cite{De96}, 1796, p.\;237).  He goes on to write down the formula
\[
\cos z = \frac{ \left[\cos \frac{z}{n} + \sin
\frac{z}{n}\sqrt{-1}\,\right]^n+ \left[\cos \frac{z}{n} - \sin
\frac{z}{n}\sqrt{-1}\,\right]^n } {2}
\]
as well as a similar formula for the sine function.  Next, de Prony
makes the following remark:
\begin{quote}
Je remarque maintenant qu'\`a mesure que~$A$ diminue et~$n$ augmente,
ces \'equations s'approchent de devenir
\begin{equation}
\label{e31}
\cos z = \frac{ \left[1+\frac{z\sqrt{-1}}{n}\right]^n +
\left[1-\frac{z\sqrt{-1}}{n}\right]^n } {2}
\end{equation}
(ibid.; labeling \eqref{e31} added)
\end{quote}
De Prony's formula~\eqref{e31} is correct only up to taking the
standard part of the right-hand side (for infinite~$n$).  Again de
Prony handles the conceptual difficulty of dealing with infinite and
infinitesimal numbers by \emph{suppressing} limits or standard parts.
Note that both of de Prony's formulas are taken verbatim from (Euler
\cite{Eu48}, 1748, \S133 -- \S138).%
\footnote{Schubring lodges the following claim concerning de\;Prony:
``The break with previous tradition, which was probably the most
visible to his contemporaries, was the exclusion and rejection of
infiniment petits by the analytic method.  In de Prony the infiniment
petits were excluded from the foundational concepts of his teaching by
simply not being mentioned; only in a heading did they appear in a
quotation, as `so-called analysis of the infinitely small
quantities'{}'' (Schubring \cite{Sc05}, 2005, p.\;289).  Schubring's
assessment of de Prony's attitude toward infinitesimals seems about as
apt as his assessment of Cauchy's; see (B\l aszczyk et
al.\;\cite{17e}, 2017).}

It is reasonable to assume that de Prony's criticism of Cauchy's
teaching of prospective engineers had to do with what Prony saw as
excessive fussiness in dealing with what came to be viewed later as
conceptual difficulties of passing to the limit, i.e., taking the
standard part.  Note that in the comment by de Prony cited at the
beginning of this section, he does \emph{not} criticize Cauchy for not
using infinitesimals, but merely for excessive emphasis on technical
detail involving \emph{limites}.  Therefore Gilain's claim to the
contrary cited at the beginning of Section~\ref{s37} amounts to
massaging the evidence by putting a tendentious spin on de Prony's
criticism.

\subsection{Foundations, limits, and infinitesimals}

Can one claim that Cauchy established the foundations of analysis on
the concept of infinitesimal?  

The notions of infinitesimal, limit, and variable quantity are all
fundamental for Cauchy.  One understands them only by the definition
which explains how they interact.  If Cauchy established such
foundations it was on the concept of a variable quantity, as analyzed
in Section~\ref{s24}.

Can one claim that Cauchy conferred upon \emph{limite} a central role
in the architecture of analysis?  The answer is affirmative if one
takes note of the frequency of the occurrence of the term in Cauchy's
oeuvre; similarly, Cauchy conferred upon infinitesimals a central role
in the said architecture.  

A more relevant issue, however, is the precise meaning of the term
\emph{limite} as used by Cauchy.  As we saw in Section~\ref{s21} he
used it in the \emph{differential} calculus in a sense closer to the
\emph{standard part function} than to any limit concept in the context
of an Archimedean continuum; and as we saw in Section~\ref{s22c}, he
used it in the \emph{integral} calculus in a sense closer to the
\emph{standard part} than any Archimedean counterpart.

Did Cauchy ever seek a justification of infinitesimals in terms of
limits?  Hardly so, since he expressed both concepts in terms of a
primitive notion of variable quantity.  In applications of analysis,
Cauchy makes no effort to justify infinitesimals in terms of limits.

\subsection{Cauchy's A-track arguments}
\label{s32}

Let us examine in more detail the issue of~$\epsilon,\delta$ arguments
in Cauchy, as found in \cite[Section~2.3, Theorem~1]{Ca21} (already
mentioned in Section~\ref{s26}).  Cauchy seeks to show that if the
difference~$f(x+1)-f(x)$ converges towards a certain limit~$k$, for
increasing values of~$x$, then the ratio~$\frac{f(x)}{x}$ converges at
the same time towards the same limit; see \cite[p.\;35]{BS}.

Cauchy chooses~$\epsilon>0$, and notes that we can give the number~$h$
a value large enough so that, when~$x$ is equal to or greater
than~$h$, the difference~$f(x+1)-f(x)$ is always contained
between~$k-\epsilon$ and \mbox{$k+\epsilon$}.  Cauchy then arrives at
the formula
\[
\frac{f(h+n)-f(h)}{n}=k+\alpha,
\]
where~$\alpha$ is a quantity contained between the limits~$-\epsilon$
and~$+\epsilon$, and eventually obtains that the
ratio~$\frac{f(x)}{x}$ has for its limit a quantity contained
between~$k-\epsilon$ and~$k+\epsilon$.

This is a fine sample of a prototype of an~$\epsilon,\delta$ proof in
Cauchy.  However, as pointed out by Sinkevich, Cauchy's proofs are all
missing the tell-tale sign of a modern proof in the tradition of the
Weierstrassian \emph{Epsilontik}, namely exhibiting an explicit
functional dependence of~$\delta$ (or in this case~$h$) on~$\epsilon$
(Sinkevich \cite{Si16}, 2016).

One of the first occurrences of a modern definition of continuity in
the style of the \emph{Epsilontik} can be found in Schwarz's summaries
of 1861 lectures by Weierstrass; see (Dugac \cite{Du73}, 1973,
p.\;64), (Yushkevich \cite{Yu86}, 1986, pp.\;74--75).  This definition
is a verbal form of a definition featuring a correct quantifier order
(involving alternations of quantifiers).

The salient point here is that this sample of Cauchy's work has no
bearing on Cauchy's infinitesimals.  Nor does it imply that
infinitesimals are merely variables tending to zero, since the term
\emph{infinitely small} does not occur in this proof at all.  Nor does
Cauchy's argument show that he thought of limits in anything
resembling post-Weierstrassian terms since his recurring definition of
limit routinely falls back on the primitive notion of a variable
quantity, rather than on any form of an alternating quantifier string,
whether verbal or not.

\subsection{Lacroix, Laplace, and Poisson}
\label{s22}

The Bradley--Sandifer edition quotes a revealing comment of Cauchy's
on the importance of infinitesimals.  The comment is found in Cauchy's
introduction:
\begin{quote}
In speaking of the continuity of functions, I could not dispense with
a treatment of the principal properties of infinitely small
quantities, properties which serve as the foundation of the
infinitesimal calculus.  (Cauchy as translated in \cite[p.\;1]{BS})
\end{quote}
Bradley and Sandifer then go on to note: ``It is interesting that
Cauchy does \emph{not} also mention limits here'' (ibid., note~6;
emphasis added).

The circumstances of the publication of the 1821 \emph{Cours
d'Analyse} indicate that attaching fundamental importance to
infinitesimals rather than limits (noted by Bradley and Sandifer) was
Cauchy's personal choice, rather than being dictated by the
constraints of his teaching at the \emph{\'Ecole Polytechnique}.
Indeed, unlike Cauchy's later textbooks, his 1821 book was not
commissioned by the \emph{\'Ecole} but was rather written upon the
personal request of Laplace and Poisson, as acknowledged in (Gilain
\cite{Gi89}, 1989, note\;139).

Sinaceur points out that Cauchy's definition of limit resembles, not
that of Weierstrass, but rather that of Lacroix%
\footnote{As a student at the \emph{Polytechnique}, Cauchy attended
Lacroix's course in analysis in 1805; see (Belhoste \cite{Be91}, 1991,
p.\;10, 243).}
dating from 1810 (see \cite[p.\;108--109]{Si73}).%
\footnote{Sinaceur explicitly denies Cauchy the honor of having
published the first arithmetic definition of limits, by writing: ``Or,
1) l'{\'e}psilonisation n'est pas l'{\oe}uvre de Cauchy, mais celle de
Weierstrass ; \ldots{} on ne peut dire qu'il en donne une
d{\'e}finition purement arithm{\'e}tique ou purement analytique.  Sa
d{\'e}finition \ldots{} n'enveloppe pas moins d'intuition
g{\'e}om{\'e}trique que celle contenue dans le \emph{Trait{\'e}} de
Lacroix\ldots{}''}
This is acknowledged in (Grabiner \cite{Gr81}, 1981, p.\;80).

Cauchy's kinematic notion of limit was expressed, like his notion of
infinitesimal~$\alpha$, in terms of a primitive notion of
\emph{variable quantity} (see Section~\ref{s24}).  Thus, Cauchy's
comment that when a variable becomes an infinitesimal~$\alpha$, the
limit of such a variable is zero, can be interpreted in two ways.  It
can be interpreted in the context of an Archimedean continuum.
Alternatively, it could be interpreted as the statement that the
assignable part of~$\alpha$ is zero, in the context of a Bernoullian
(i.e., infinitesimal-enriched) continuum, or in modern terminology,
that the \emph{standard part} of~$\alpha$ is zero; see
Section~\ref{s71}.

\section{Modern infinitesimals in relation to Cauchy's procedures}
\label{s71}

While set-theoretic justifications for either A-track or B-track
modern framework are obviously not to be found in Cauchy, Cauchy's
\emph{procedures} exploiting infinitesimals find closer proxies in
Robinson's framework for analysis with infinitesimals than in a
Weierstrassian framework.  In this section we outline a set-theoretic
construction of a hyperreal extension~$\R\hookrightarrow\astr$, and
point out specific similarities between procedures using the
hyperreals, on the one hand, with Cauchy's procedures, on the other.

Let~$\R^{\N}$ denote the ring of sequences of real numbers, with
arithmetic operations defined termwise.  Then we have
$\astr=\R^{\N}\!/\,\text{MAX}$ where MAX is the maximal ideal
consisting of all ``negligible'' sequences~$(u_n)$.  Here a sequence
is negligible if it vanishes for a set of indices of full
measure~$\xi$, namely,~$\xi\big(\{n\in\N\colon u_n=0\}\big)=1$.  Here
$\xi\colon \mathcal{P}(\N)\to \{0,1\}$ is a finitely additive
probability measure taking the value~$1$ on cofinite sets,
where~$\mathcal{P}(\N)$ is the set of subsets of~$\N$.  The
subset~$\mathcal{F}_\xi\subseteq\mathcal{P}(\N)$ consisting of sets of
full measure~$\xi$ is called a nonprincipal ultrafilter.  These
originate with (Tarski \cite{Ta30}, 1930).  The set-theoretic
presentation of a Bernoullian continuum (see Section~\ref{s12})
outlined here was therefore not available prior to that date.

The field~$\R$ is embedded in~$\astr$ by means of constant sequences.
The subring $\hr\subseteq\astr$ consisting of the finite elements
of~$\astr$ admits a map~$\st$ to~$\R$, known as \emph{standard part}
\begin{equation}
\label{e61}
\st\colon \hr\to\R,
\end{equation}
which rounds off each finite hyperreal number to its nearest real
number (the existence of such a map \st\;is the content of the
\emph{standard part principle}).  This enables one, for instance, to
define the derivative of $s=f(t)$ as
\begin{equation}
\label{e62}
f'(t)=\frac{ds}{dt}= \st\left(\frac{\Delta s}{\Delta t}\right)
\end{equation}
(here~$\Delta s\ne0$ is infinitesimal) which parallels Cauchy's
definition of derivative (see equation~\eqref{e21} in
Section~\ref{s21}) more closely than any \emph{Epsilontik} definition.
Limit is similarly defined in terms of \st{}, e.g., by setting
\begin{equation}
\label{e83}
\lim_{t\to0}f(t)=\st(f(\epsilon))
\end{equation}
where~$\epsilon$ is a nonzero infinitesimal, in analogy with Cauchy's
limit as analyzed in Section~\ref{s12c}.  For additional details on
Robinson's framework see e.g., \cite{17f}.

\section{Conclusion}

The oft-repeated claim (as documented e.g., in \cite{17a}; \cite{18e})
that ``Cauchy's infinitesimal is a variable with limit~$0$'' (see
Gilain's comment cited in Section~\ref{s34}) is a reductionist view of
Cauchy's foundational stance, at odds with much compelling evidence in
Cauchy's writings, as we argued in Sections~\ref{s2} and \ref{s5}.

Gilain, Siegmund-Schultze, and some other historians tend to adopt a
butterfly model for the development of analysis, to seek proxies for
Cauchy's procedures in a default modern Archimedean framework, and to
view his infinitesimal techniques as an evolutionary dead-end in the
history of analysis.  Such an attitude was criticized by
Grattan-Guinness, as discussed in Section~\ref{s1b}.  The fact is
that, while Cauchy did use an occasional epsilon in an Archimedean
sense, his techniques relying on infinitesimals find better proxies in
a modern framework exploiting a Bernoullian continuum.

Robinson first proposed an interpretation of Cauchy's
\emph{procedures} in the framework of a modern theory of
infinitesimals in \cite{Ro66} (see Section~\ref{s1}).  A
\emph{set-theoretic foundation} for infinitesimals could not have been
provided by Cauchy for obvious reasons, but Cauchy's \emph{procedures}
find closer proxies in modern infinitesimal frameworks than in modern
Archimedean ones.

\section*{Acknowledgments} 

We are grateful to Peter Fletcher for helpful suggestions.  We thank
Reinhard Siegmund-Schultze for bringing his review \cite{Si09} to our
attention.

\end{document}